\documentclass[preprint,3p,12pt,pdf]{elsarticle}

\usepackage{lineno,hyperref}
\modulolinenumbers[5]

\usepackage{hyperref}

\usepackage{epstopdf}

\usepackage{tikz}
\usetikzlibrary{arrows}

\usepackage{pst-node}
\usepackage{tikz-cd}

\usepackage[shellescape]{gmp}

\usepackage{mathrsfs}
\usepackage{amssymb}
\usepackage{amsfonts}
\usepackage{latexsym}
\usepackage{mathtools}
\usepackage{xcolor}
\usepackage{wrapfig}
\usepackage{floatflt}

\usepackage{mathtools}
\usepackage{extarrows}

\usepackage{graphicx}

\usepackage{subcaption}

\def\lb{\label}

\newcommand{\er}[1]{\textrm{(\ref{#1})}}

\def\a{\alpha}

\def\G{\Gamma}

   \def\vt{\vartheta}        

    \def\R{{\mathbb R}}       
    \def\N{{\mathbb N}}   

\let\ge\geqslant                 

\def\iy{\infty}

\def\el2{\ell^{\,2}}             \def\1{1\!\!1}

\def\Im{\mathop{\mathrm{Im}}\nolimits}

\let\ge\geqslant

\newcommand{\ca}{\begin{cases}}
	\newcommand{\ac}{\end{cases}}
\newcommand{\ma}{\begin{pmatrix}}
	\newcommand{\am}{\end{pmatrix}}
\def\eq{\begin{equation}}
	\def\qe{\end{equation}}
\def\[{\begin{equation}}
	\def\]{\end{equation}}

\begin{document}
	
	\begin{frontmatter}

		\title{A formula for the periodic multiplier in left tail asymptotics for supercritical branching processes in the Schr\"oder case}

		\date{\today}

		\author
		{Anton A. Kutsenko}
	
		\address{Mathematical Institute for Machine Learning and Data
Science, KU Eichst\"att--Ingolstadt, Germany; email: akucenko@gmail.com}
	
	\begin{abstract}
		It is known that the left tail asymptotic for supercritical branching processes in the Schr\"oder case satisfies a power law multiplied by some multiplicatively periodic function. We provide an explicit expression for this periodic function.  
	\end{abstract}


	
\end{frontmatter}


{\section{Introduction}\lb{sec0}}

We consider a simple Galton-Watson branching process $Z_t$ in the supercritical case with the minimum family size $1$ - the so-called Schr\"oder case. Denoting the mean of offspring distribution by $E$, one may define the {\it martingal limit} $W=\lim_{t\to+\iy}E^{-t}Z_t$. It is proven in \cite{BB1} that $p(x)$, the density of $W$, has an asymptotic
\[\lb{000}
 p(x)=x^{\a}V(x)+o(x^{\a}),\ \ \ x\to0,
\]
with some explicit $\a$ and a continuous, positive, multiplicatively periodic function, $V$, with period $E$. Further references to this asymptotic are always based on the principal work \cite{BB1} and do not provide any formula for $V(x)$. We will try to fill this gap. Our explicit expression for $V(x)$ is written in terms of Fourier coefficients of $1$-periodic Karlin-McGregor function $K^*(z)$, see \er{015}, and some values of $\G$-function, see \er{310}. The main result is given in \er{a020}. 

{\section{Main results}\lb{sec1}}

Let us recall some known facts about branching processes, including some recent results obtained in \cite{K}. The Galton--Watson process is defined by
\[\lb{001}
 X_{t+1}=\sum_{j=1}^{X_t}\xi_{j,t},\ \ \ X_0=1,\ \ \ t\in\N\cup\{0\},
\]
where all $\xi_{j,t}$ are independent and identically-distributed natural number-valued random variables with the probability generating function
\[\lb{002}
 P(z):=\mathbb{E}z^{\xi}=p_0+p_1z+p_2z^2+p_3z^3+....
\]
It is well known, see \cite{H}, that 
\[\lb{003}
 \mathbb{E}z^{X_t}=\underbrace{P\circ...\circ P}_{t}(z)=p_{t0}+p_{t1}z+p_{t2}z^2+...,\ \ \ t\in\N.
\]
For simplicity, let us assume that $P$ is a polynomial of degree $N\ge2$, which is commonly the case in practice. Another assumption is $p_0=0$ and $p_1\ne0$. To avoid complications due to periodicity, we assume also the existence of $k$ such that $p_kp_{k+1}\ne0$. This last assumption is not important for the {\it martingal limit} but, for simplicity, we would like to keep it because our derivation of the main asymptotic is based on using another limit distribution for which the assumption is significant. It is useful to note that the case $p_0\ne0$ can usually be reduced to the (supercritical) case $p_0=0$ with the help of Harris-Sevastyanov transformation, see \cite{B1}. Under the assumptions $p_0=0$ and $0<p_1<1$, we can define
\[\lb{004}
 \Phi(z)=\lim_{t\to\iy}p_1^{-t}\underbrace{P\circ...\circ P}_{t}(z),
\]
which is analytic at least for $|z|<1$. The function $\Phi$ satisfies the Scr\"oder-type functional equation
\[\lb{005}
 \Phi(P(z))=p_1\Phi(z),\ \ \ \Phi(0)=0,\ \ \ \Phi'(0)=1.
\] 
Due to \er{002}-\er{004}, the Taylor coefficients $\varphi_n$ of
\[\lb{006}
 \Phi(z)=z+\varphi_2z^2+\varphi_3z^3+...
\]
describe the relative limit densities of the number of descendants 
\[\lb{007}
 \varphi_n=\lim_{t\to+\iy}\frac{p_{tn}}{p_{t1}}=\lim_{t\to+\iy}\frac{\mathbb{P}\{Z_t=n\}}{\mathbb{P}\{Z_t=1\}}.
\]
In \cite{K}, it is shown that the densities have the asymptotic
\[\lb{0070}
 \varphi_n\simeq n^{\frac{-\ln (Ep_1)}{\ln E}}K_0(\frac{-\ln n}{\ln E})+o(n^{\frac{-\ln (Ep_1)}{\ln E}}),
\]
where the first moment
\[\lb{008}
 E:=P'(1)=\sum_{j=1}^{N}jp_j>1,
\]
and $1$-periodic function $K_0(z)$ is given by
\[\lb{310}
 K_0(z)=\sum_{m=-\iy}^{+\iy}\frac{\theta_me^{2\pi\mathbf{i}mz}}{\Gamma(-\frac{2\pi\mathbf{i}m+\ln p_1}{\ln E})}.
\]
The numbers $\theta_m$ are Fourier coefficients of the $1$-periodic function
\[\lb{015}
 K^*(z):=\Phi(\Pi(E^z))p_1^{-z}=\sum_{m=-\iy}^{+\iy}\theta_me^{2\pi \mathbf{i} mz},
\] 
defined on the strip $\{z:\ |\Im z|<\vt^*\}$ with some $\vt^*>0$. The width $2\vt^*$ of the strip depends on the geometric properties of the Julia set for the polynomial $P$, see \cite{K}. This width affects the rate of exponential decaying of the coefficients $\theta_m$. The function $\Pi(z)$ is given by
\[\lb{013}
 \Pi(z):=\lim_{t\to+\iy}\underbrace{P\circ...\circ P}_{t}(1-E^{-t}z).
\]
This function is entire satisfying the Poincar\'e-type functional equation
\[\lb{014}
 P(\Pi(z))=\Pi(Ez),\ \ \ \Pi(0)=1,\ \ \ \Pi'(0)=-1.
\]
In addition to \er{007}, there is another scaling of the coefficients $p_{tn}$, namely
\[\lb{a015}
 p(x)=\lim_{t\to+\iy}E^tp_{t,[xE^t]},\ \ \ x\in\R,
\]
where $[\cdot]$ denotes the integer part. It is assumed that $p(x)=0$ for $x<0$. On the one side, we have the Riemann sum
\[\lb{a016}
 \underbrace{P\circ...\circ P}_{t}(e^{\frac{-\mathbf{i}y}{E^t}})=E^{-t}\sum_{n=0}^{+\iy} E^{t}p_{tn}e^{\frac{-n\mathbf{i}y}{E^t}}\to \int_{-\iy}^{+\iy}p(x)e^{-\mathbf{i}yx}dx.
\]
On the other side, the well known expansion for the exponent gives
\[\lb{a017}
 \underbrace{P\circ...\circ P}_{t}(e^{\frac{-\mathbf{i}y}{E^t}})\sim \underbrace{P\circ...\circ P}_{t}(1-{\frac{\mathbf{i}y}{E^t}})\to\Pi(\mathbf{i}y).
\]
Thus, combining \er{a016}, \er{a017}, and applying the inverse Fourier transform, we obtain
\[\lb{a018}
 p(x)=\frac1{2\pi}\int_{-\iy}^{+\iy}\Pi(\mathbf{i}y)e^{\mathbf{i}yx}dy.
\]
Using $1$-periodicity of $K_0$, identity $p_{1t}=p_1^t$, see \er{004} and \er{005}, along with  \er{007} and \er{0070}, we obtain
\[\lb{a019}
  E^tp_{t,[xE^t]}=E^tp_1^t\frac{p_{t,[xE^t]}}{p_{1t}}\sim (Ep_1)^t(xE^t)^{\frac{-\ln Ep_1}{\ln E}}   K_0(\frac{-\ln xE^t}{\ln E})=x^{\frac{-\ln Ep_1}{\ln E}}K_0(\frac{-\ln x}{\ln E}).
\]
Finally, \er{a019}, \er{a015}, and \er{000} gives 
\[\lb{a020}
 V(x)= K_0(\frac{-\ln x}{\ln E}).
\]
A computation of $K_0$ can be based on efficient numerical schemes presented in \cite{K}. Fourier transform \er{a018}, as well as Fourier coefficients $\theta_m$, see \er{015}, can be computed efficiently with the help of FFT. We have prepared two numerical examples showing the comparison of exact values of the normalized density and the first asymptotic term, see Fig. \ref{fig1}.
\begin{figure}[h]
	\center{\includegraphics[width=0.9\linewidth]{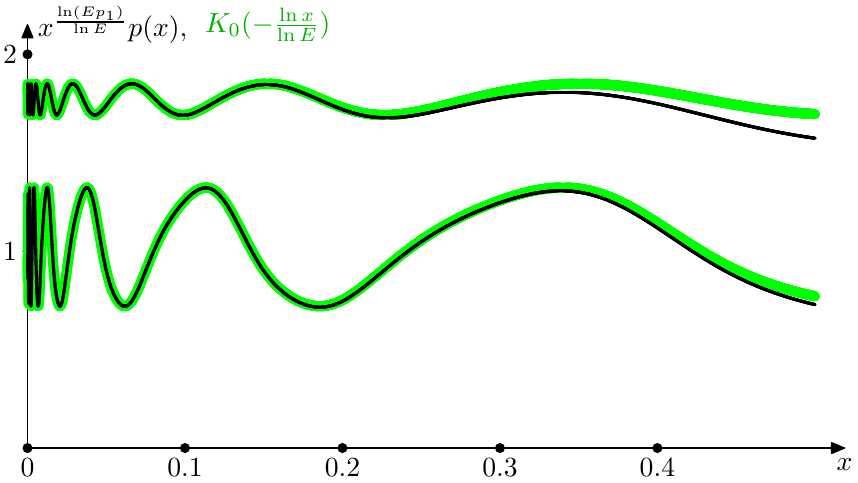}}
	\caption{Two examples: the cases $p_1=0.1$, $p_2=0.5$, $p_3=0.4$ (upper curves), and $p_1=0.1$, $p_2=0.1$, $p_3=0.5$, $p_4=0.3$.}\lb{fig1}
\end{figure} 

\section*{Acknowledgements} 
This paper is a contribution to the project M3 of the Collaborative Research Centre TRR 181 "Energy Transfer in Atmosphere and Ocean" funded by the Deutsche Forschungsgemeinschaft (DFG, German Research Foundation) - Projektnummer 274762653. 

\section*{Data availability statement}

Data sharing not applicable to this article as no datasets were generated or analyzed during the current study. The program code for procedures and functions considered in the article is available from the corresponding author on reasonable request.

\end{document}